Р. Ибрагимов, Ш. Шарахметов

# МОМЕНТНЫЕ НЕРАВЕНСТВА
# ДЛЯ СИММЕТРИЧЕСКИХ СТАТИСТИК

В работах [1], [2] были получены аналоги моментных неравенств Хинчина и Розенталя для симметрических статистик вида

$$T_n = \sum_{1 \le i < j \le n} Y(X_i, X_j) \ .$$

Цель настоящей статьи - обобщить эти результаты на случай симметрических статистик произвольного порядка. Кроме аналогов неравенств Хинчина и Розенталя, в работе доказаны двусторонние оценки для моментов неотрицательных статистик.

Пусть, так же, как в [1], [2], $X_1,...,X_n$ - независимые одинаково распределенные случайные величины (с.в.), принимающие значения в некотором измеримом пространстве $(\Im, A)$, $1 \le m \le n$, $p \ge 1$.

Введем в рассмотрение класс $F(p)$ функций $Y: \Im^m \to R$, удовлетворяющих условиям

$$Y(X_1,X_2,...,X_m) = Y(X_{\pi(1)}, X_{\pi(2)}, ..., X_{\pi(m)}) \quad (\text{п.н.}) \tag{1}$$

для всех перестановок $\pi$ множества $\{1,2,...,m\}$,

$$E|Y(X_1,X_2,...,X_m)|^p < \infty \ . \tag{2}$$

Через $G(p)$ будем обозначать подмножество $F(p)$, состоящее из функций $Y$, для которых

$$E(Y(X_1,X_2,...,X_m)/X_1,...,X_{m-1}) = 0 \quad (\text{п.н.}). \tag{3}$$

Для $Y \in F(p)$ положим

$$T_n = \sum_{1 \le i_1 < ... < i_m \le n} Y(X_{i_1}, X_{i_2}, ..., X_{i_m}) \ .$$

Всюду в дальнейшем $A(p)$, $B(p)$ обозначают константы, зависящие только от $p$, $A(p,m)$, $B(p,m)$ - константы, зависящие только от $p$ и $m$, вообще говоря, разные в разных местах.

**Теорема 1.** Если $Y \in F(p)$ - неотрицательная функция, $p \ge 1$, то

$$A(p,m) \max_{k=0,m} n^{p(m-k)+k} E(E(Y(X_1,X_2,...,X_m)/X_1,...,X_k))^p \le ET_n^p \le$$

$$\leq B(p,m) \max_{k=0,m} n^{p(m-k)+k} E(E(Y(X_1,X_2,...,X_m)/X_1,...,X_k))^p \quad (4)$$

(здесь и всюду далее $E(\cdot/X_1,...,X_k)=E(\cdot)$ при $k=0$).

**Теорема 2.** Если $p \geq 2$, $Y \in G(p)$, то

$$A(p,m) \max_{k=0,m} n^{p/2(m-k)+k} E(E(Y^2(X_1,X_2,...,X_m)/X_1,...,X_k))^{p/2} \leq E|T_n|^p \leq$$

$$\leq B(p,m) \max_{k=0,m} n^{p/2(m-k)+k} E(E(Y^2(X_1,X_2,...,X_m)/X_1,...,X_k))^{p/2}. \quad (5)$$

**Теорема 3.** Если $p \geq 2$, $Y \in G(p)$, то

$$A(p,m) E\left( \sum_{1 \leq i_1 < ... < i_m \leq n} Y^2(X_{i_1}, X_{i_2},...,X_{i_m}) \right)^{p/2} \leq E|T_n|^p \leq$$

$$\leq B(p,m) E\left( \sum_{1 \leq i_1 < ... < i_m \leq n} Y^2(X_{i_1}, X_{i_2},...,X_{i_m}) \right)^{p/2}. \quad (6)$$

**Замечание.** Отметим, что, вообще говоря, каждый член, участвующий в выражении

$$\varphi_n = \max_{k=0,m} n^{p/2(m-k)+k} E(E(Y^2(X_1,X_2,...,X_m)/X_1,...,X_k))^{p/2}, \; p>2,$$

(теорема 2) и, следовательно, в выражении

$$\psi_n = \max_{k=0,m} n^{p(m-k)+k} E(E(Y(X_1,X_2,...,X_m)/X_1,...,X_k))^p, \; p>1,$$

(теорема 1) является существенным. Действительно, пусть $0 \leq k \leq m$, $a=1/p-n^{-p/2}$, с.в. $X_1,...,X_m$ имеют показательное распределение с параметром $1$,

$$Y(x_1,...,x_m) = \sum_{\pi} \left[ exp\left( a\sum_{i=1}^{k} x_{\pi(i)} - \sum_{i=k+1}^{m} x_{\pi(i)} \right) + \right.$$

$$\sum_{l=1}^{m} (-1)^l \sum_{1 \leq j_1 < ... < j_l \leq m} \prod_{s=1}^{l} (1/(1-a \cdot sign(k+0.5-j_s))) exp\left( a \sum_{\substack{i=1 \\ i \neq j_1,...,j_l}}^{k} x_{\pi(i)} - \right.$$

$$a \sum_{\substack{i=k+1 \\ i \neq j_1,...,j_l}}^{m} x_{\pi(i)})],$$

где внешнее суммирование производится по всем перестановкам $\pi$ множества $\{1,2,...,m\}$. Нетрудно видеть, что $Y \in G(p)$. При $n \to \infty$ выражения

$$n^{p/2(m-l)+l} E(E(Y^2(X_1,X_2,...,X_m)/X_1,...,X_l))^{p/2}, \; l=0,...,k,$$

имеют порядок роста $n^{pm/2+l}$, а выражения

$$n^{p/2(m-l)+l} E(E(Y^2(X_1,X_2,...,X_m)/X_1,...,X_l))^{p/2}, \; l=k+1,...,m,$$

имеют порядок роста $n^{p/2(m+k)-(p/2-1)l}$. Следовательно, порядок $\varphi_n$ определяется членом $n^{p/2(m-k)+k} E(E(Y^2(X_1,X_2,...,X_m)/X_1,...,X_k))^{p/2}$, что доказывает его существенность.

Для доказательства теорем нам понадобятся следующие леммы.

**Лемма 1.** [3, теорема 5.1]. Пусть $1 \leq p < \infty$, $\Im_0 = (\Omega, \emptyset) \subseteq \Im_1 \subseteq \Im_2 \subseteq ... \subseteq \Im_n$ - возрастающая последовательность $\sigma$-алгебр на некотором вероятностном пространстве $(\Omega, \Im, P)$, $X_k, k=1,...,n$ - последовательность неотрицательных $\Im_k$-измеримых случайных величин таких, что $EX_k^p < \infty$. Тогда

$$E(\sum_{k=1}^{n} X_k)^p \leq B(p) \max(\sum_{k=1}^{n} EX_k^p, \; E(\sum_{k=1}^{n} E(X_k/\Im_{k-1}))^p).$$

**Лемма 2.** Пусть $\{X_n\}_{n=1}^{\infty}$ - последовательность неотрицательных с.в. со значениями в некотором измеримом пространстве $(\Im, A)$,

$$A_p = \sum_{n=1}^{\infty} EX_n^p < \infty,$$

$$B_\gamma = (\sum_{n=1}^{\infty} EX_n^\gamma)^{1/\gamma} < \infty.$$

Тогда

1) $A_s \leq (A_p^{s-\gamma} B_\gamma^{\gamma(p-s)})^{1/(p-\gamma)}$ для всех $1 \leq \gamma < s < p$,

2) $A_p B_\gamma^s \leq max(A_{p+s}, B_\gamma^{p+s})$ для всех $p > \gamma$, $s \geq 0$,

3) $max(A_s, B_\gamma^s) \leq (max(A_p, B_\gamma^p))^{s/p}$ для всех $1 \leq \gamma < s < p$.

Доказательство аналогично доказательству леммы в [2].

Из леммы 2 вытекает, что для всех неотрицательных функций $Y \in F(2s)$ при $k, l = 0, ..., m$, $s \geq 1$ имеют место неравенства:

$$n^{k+l+s(2m-k-l)} E(E(Y(X_1,...,X_m)/X_1,...,X_k))^s E(E(Y(X_1,...,X_m)/X_1,...,X_l))^s \leq$$

$$\leq B(s,m) max(n^{k+2s(m-k)} E(E(Y(X_1,...,X_m)/X_1,...,X_k))^{2s}, n^{l+2s(m-l)} E(E(Y(X_1,...,X_m)/X_1,...,X_l))^{2s},$$

$$n^{2sm}(EY(X_1,...,X_m)^{2s}). \tag{7}$$

<u>Доказательство теоремы 1.</u> Пусть $Y_{i_1,i_2,...,i_m} = Y(X_{i_1}, X_{i_2}, ..., X_{i_m})$, $1 \leq i_1 < ... < i_m \leq n$. Используя неотрицательность функций $Y$ и неравенство Иенсена, имеем, что

$$ET_n^p \geq \max_{k=1,m-1} \left( \sum_{1 \leq i_1 < ... < i_m \leq n} EY_{i_1,i_2,...,i_m}^p, \left( \sum_{1 \leq i_1 < ... < i_m \leq n} EY_{i_1,i_2,...,i_m} \right)^p, \right.$$

$$\sum_{k+1 \leq i_{k+1} < ... < i_m \leq n} E\left( \sum_{1 \leq i_1 < ... < i_k \leq i_{k+1}-1} Y_{i_1,i_2,...,i_m} \right)^p \right) \geq$$

$$\geq \max_{k=1,m-1} \left( \sum_{1 \leq i_1 < ... < i_m \leq n} EY_{i_1,i_2,...,i_m}^p, \left( \sum_{1 \leq i_1 < ... < i_m \leq n} EY_{i_1,i_2,...,i_m} \right)^p, \right.$$

$$\sum_{k+1 \leq i_{k+1} < ... < i_m \leq n} E\left( \sum_{1 \leq i_1 < ... < i_k \leq i_{k+1}-1} E(Y_{i_1,i_2,...,i_m}/X_{i_{k+1}},...,X_{i_m}) \right)^p \right) \geq$$

$$\geq A(p,m) \max_{k=0,m} n^{p(m-k)+k} E(E(Y(X_1,X_2,...,X_m)/X_1,...,X_k))^p.$$

Левое неравенство (4) доказано.

Положим последовательно

$$Y_{i_k,...,i_m} = \sum_{i_{k-1}=k-1}^{i_k-1} Y_{i_{k-1},i_k,...,i_m}, \quad k \leq i_k < ... < i_m \leq n, \quad k=2,...,m.$$

Очевидно, что с.в. $Y_{i_m}$, $i_m=m,...,n$ измеримы относительно $\sigma$-алгебр $\sigma(X_1,X_2,...,X_{i_m})$, с.в. $Y_{i_k,...,i_m}$, $i_k=k,...,i_{k+1}-1$ измеримы относительно $\sigma$-алгебр $\sigma(X_1,X_2,...,X_{i_k},X_{i_{k+1}},...,X_{i_m})$, $k+1 \leq i_{k+1} < ... < i_m \leq n$, $k=1,...,m-1$.

Применяя лемму 1, получаем, что

$$ET_n^p \leq B(p) \max\left( \sum_{i_m=m}^{n} EY_{i_m}^p, E\left( \sum_{i_m=m}^{n} E(Y_{i_m}/X_1,...,X_{i_m-1})\right)^p \right), \qquad (8)$$

$$EY_{i_k,...,i_m}^p \leq B(p) \max\left( \sum_{i_{k-1}=k-1}^{i_k-1} EY_{i_{k-1},...,i_m}^p, E\left( \sum_{i_{k-1}=k-1}^{i_k-1} E(Y_{i_{k-1},...,i_m}/X_1,...,X_{i_{k-1}-1},X_{i_k},...,X_{i_m})\right)^p \right), \qquad (9)$$

$k \leq i_k < ... < i_m \leq n$, $k=2,...,m$.

Из (8) и (9) вытекает, что

$$ET_n^p \leq B(p) \max_{k=1,m-1} \Bigg( \sum_{1 \leq i_1 < ... < i_m \leq n} EY_{i_1,i_2,...,i_m}^p, E\left( \sum_{1 \leq i_1 < ... < i_m \leq n} E(Y_{i_1,i_2,...,i_m}/X_{i_1},X_{i_2},...,X_{i_{m-1}})\right)^p,$$

$$\sum_{k+1 \leq i_{k+1} < ... < i_m \leq n} E\left( \sum_{1 \leq i_1 < ... < i_k \leq i_{k+1}-1} E(Y_{i_1,i_2,...,i_m}/X_{i_1},X_{i_2},...,X_{i_{k-1}},X_{i_{k+1}},...,X_{i_m})\right)^p \Bigg) \leq$$

$$\leq B(p,m) \max_{k=1,m-1} \Bigg( n^m EY_{1,2,...,m}^p, n^{p-1} \sum_{i_m=m}^{n} E\left( \sum_{1 \leq i_1 < ... < i_{m-1} \leq i_m-1} E(Y_{i_1,i_2,...,i_m}/X_{i_1},X_{i_2},...,X_{i_{m-1}})\right)^p,$$

$$\sum_{k+1 \leq i_{k+1} < ... < i_m \leq n} i_{k+1}^{p-1} \sum_{i_k=k}^{i_{k+1}-1} E\left( \sum_{1 \leq i_1 < ... < i_{k-1} \leq i_k-1} E(Y_{i_1,i_2,...,i_m}/X_{i_1},X_{i_2},...,X_{i_{k-1}},X_{i_{k+1}},...,X_{i_m})\right)^p \Bigg). \qquad (10)$$

Предполагая правое неравенство (4) доказанным при $m \leq l$, из (10) легко получаем, что оно справедливо при $m=l+1$. Поскольку при $m=1$ (4) имеет место, по принципу индукции утверждение теоремы полностью доказано.

<u>Доказательство теоремы 2.</u> Пусть с.в. $Y_{i_k,...,i_m}$, $k \leq i_k < ... < i_m \leq n$, $k=1,...,m$ те же, что и в доказательстве теоремы 1. Легко видеть, что при условиях теоремы 2 с.в. $Y_{i_m}$, $i_m=m,...,n$ образуют мартингал-разность относительно $\sigma$-алгебр $\sigma(X_1,X_2,...,X_{i_m})$, с.в. $Y_{i_k,...,i_m}$, $i_k=k,...,i_{k+1}-1$ образуют мартингал-разность относительно $\sigma$-алгебр $\sigma(X_1,X_2,...,X_{i_k},X_{i_{k+1}},...,X_{i_m})$, $k+1 \leq i_{k+1} < ... < i_m \leq n$, $k=1,...,m-1$.

Используя аналог неравенства Розенталя для мартингалов [4, теорема 2.1.5], имеем, что

$$A(p) \max\left( \sum_{i_m=m}^{n} E|Y_{i_m}|^p, E\left( \sum_{i_m=m}^{n} E(Y_{i_m}^2/X_1,...,X_{i_m-1})\right)^{p/2} \right) \leq E|T_n|^p \leq$$

$$\leq B(p)\, max\left(\sum_{i_m=m}^{n} E|Yi_m|^p,\ E\left(\sum_{i_m=m}^{n} E(Y i_m^2 / X_1,..., X_{i_m-1})\right)^{p/2}\right), \tag{11}$$

$$A(p) max\left(\sum_{i_{k-1}=k-1}^{i_k-1} E|Y i_{k-1},...,i_m|^p,\ E\left(\sum_{i_{k-1}=k-1}^{i_k-1} E(Y i_{k-1},...,i_m^2 / X_1,..., X_{i_{k-1}-1}, X_i,..., X_{i_m})\right)^{p/2}\right) \leq E|Y i_k,...,i_m|^p \leq$$

$$\leq B(p)\, max\left(\sum_{i_{k-1}=k-1}^{i_k-1} E|Y i_{k-1},...,i_m|^p,\ E\left(\sum_{i_{k-1}=k-1}^{i_k-1} E(Y i_{k-1},...,i_m^2 / X_1,..., X_{i_{k-1}-1}, X_{i_k},..., X_{i_m})\right)^{p/2}\right), \tag{12}$$

$k \leq i_k < ... < i_m \leq n$, $k = 2,...,m$.

Из (11), (12) вытекает, что

$$A(p)\, \max_{k=1,m-1}\left(\sum_{1\leq i_1<...<i_m\leq n} E|Y i_1,i_2,...,i_m|^p,\ E\left(\sum_{i_m=m}^{n} E(Y i_m^2 / X_1,...,X_{i_m-1})\right)^{p/2},\right.$$

$$\left.\sum_{k+1\leq i_{k+1}<...<i_m\leq n} E\left(\sum_{i_k=k}^{i_{k+1}-1} E(Y i_k,...,i_m^2 / X_1, X_2,..., X_{i_k-1}, X_{i_{k+1}},..., X_{i_m})\right)^{p/2}\right) \leq E|T_n|^p \leq$$

$$\leq B(p)\, \max_{k=1,m-1}\left(\sum_{1\leq i_1<...<i_m\leq m} E|Y i_1,i_2,...,i_m|^p,\ E\left(\sum_{i_m=m}^{n} E(Y i_m^2 / X_1,...,X_{i_m-1})\right)^{p/2},\right.$$

$$\left.\sum_{k+1\leq i_{k+1}<...<i_m\leq n} E\left(\sum_{i_k=k}^{i_{k+1}-1} E(Y i_k,...,i_m^2 / X_1, X_2,..., X_{i_k-1}, X_{i_{k+1}},..., X_{i_m})\right)^{p/2}\right). \tag{13}$$

По неравенству Иенсена

$$E\left(\sum_{i_m=m}^{n} E(Y i_m^2 / X_1,..., X_{i_m-1})\right)^{p/2} \geq \left(\sum_{i_m=m}^{n} E\, Y i_m^2\right)^{p/2} =$$

$$= \left(\sum_{1\leq i_1<...<i_m\leq n} E Y i_1,i_2,...,i_m^2\right) \geq A(p,m)\, n^{p/2m}(E Y_{1,...,m}^2)^{p/2}. \tag{14}$$

$$\sum_{k+1\leq i_{k+1}<...<i_m\leq n} E\left(\sum_{i_k=k}^{i_{k+1}-1} E(Y i_k,...,i_m^2 / X_1, X_2,..., X_{i_k-1}, X_{i_{k+1}},..., X_{i_m})\right)^{p/2} \geq$$

$$\geq \sum_{k+1\leq i_{k+1}<...<i_m\leq n} E\left(\sum_{i_k=k}^{i_{k+1}-1} E(Y i_k,...,i_m^2 / X_{i_{k+1}},..., X_{i_m})\right)^{p/2} =$$

$$= \sum_{k+1 \leq i_{k+1} < ... < i_m \leq n} E\left(\sum_{1 \leq i_1 < ... < i_k \leq i_{k+1}-1} E(Y_{i_1,...,i_m}/X^2_{i_{k+1}},...,X_{i_m})\right)^{p/2} \geq$$

$$\geq A(p,m)\, n^{kp/2+(m-k)} E(E(Y_{1,...,m}/X^2_{k+1},...,X_m))^{p/2}, \quad k=1,...,m-1. \tag{15}$$

Из (13)-(15) следует левая часть неравенства (5).

Докажем правую часть (5). Из (13) имеем, что

$$E|T_n|^p \leq B(p,m) \max_{k=1,m-1} (n^m E|Y_{1,2,...,m}|^p,$$

$$n^{p/2-1} \sum_{i_m=m}^{n} E\left(\sum_{1 \leq i_1 < ... < i_{m-1} \leq i_m-1} E(Y_{i_1,i_2,...,i_m}/X^2_{i_1},X_{i_2},...,X_{i_{m-1}})\right)^{p/2},$$

$$\sum_{k+1 \leq i_{k+1} < ... < i_m \leq n} i_{k+1}^{p/2-1} \sum_{i_k=k}^{i_{k+1}-1} E\left(\sum_{1 \leq i_1 < ... < i_{k-1} \leq i_k-1} E(Y_{i_1,i_2,...,i_m}/X^2_{i_1},X_{i_2},...,X_{i_{k-1}},X_{i_{k+1}},...,X_{i_m})\right)^{p/2} +$$

$$+ B(p,m) \max_{k=2,m-1} (n^{p/2-1} \sum_{i_m=m}^{n} E|W_{2m-2}|^{p/2}, \sum_{k+1 \leq i_{k+1} < ... < i_m \leq n} i_{k+1}^{p/2-1} \sum_{i_k=k}^{i_{k+1}-1} E|W_{2k-2}|^{p/2}), \tag{16}$$

где

$$W_{2m-2} = \sum_{s=1}^{m-1} \sum_{\substack{1 \leq i_1 < ... < i_{m-1} \leq i_m-1 \\ 1 \leq j_1 < ... < j_{m-1} \leq i_m-1 \\ i_l = j_l,\ l=1,...s-1 \\ i_s < j_s}} E(Y_{i_1,...,i_{m-1},i_m} Y_{j_1,...,j_{m-1},i_m}/X_{i_1},...,X_{i_{m-1}},X_{j_1},...,X_{j_{m-1}}),$$

$$W_{2k-2} = \sum_{s=1}^{k-1} \sum_{\substack{1 \leq i_1 < ... < i_{k-1} \leq i_k-1 \\ 1 \leq j_1 < ... < j_{k-1} \leq i_k-1 \\ i_l = j_l,\ l=1,...s-1 \\ i_s < j_s}} E(Y_{i_1,...,i_{k-1},i_{k+1},...,i_m} Y_{j_1,...,j_{k-1},i_{k+1},...,i_m}/X_{i_1},...,X_{i_{k-1}},X_{j_1},...,X_{j_{k-1}},X_{i_{k+1}},...,X_{i_m}).$$

Из теоремы 1 следует, что

$$E\left(\sum_{1 \leq i_1 < ... < i_{m-1} \leq i_m-1} E(Y_{i_1,i_2,...,i_m}/X^2_{i_1},X_{i_2},...,X_{i_{m-1}})\right)^{p/2} \leq$$

$$\leq B(p,m) \max_{s=0,m-1} i_m^{s+p/2(m-1-s)} E(E(Y_{1,...,m}^2/X_1,...,X_s))^{p/2}, \qquad (17)$$

$$E(\sum_{1\leq i_1<...<i_{k-1}\leq i_k-1} E(Y_{i_1,i_2,...,i_m}^2/X_{i_1},X_{i_2},...,X_{i_{k-1}},X_{i_{k+1}},...,X_{i_m}))^{p/2} \leq$$

$$\leq B(p,m) \max_{s=0,k-1} i_k^{s+p/2(m-1-s)} E(E(Y_{1,...,m}^2/X_1,...,X_s))^{p/2}. \qquad (18)$$

Используя неравенство Иенсена, получаем, что при $1 \leq p < 2$, $Y \in G(p)$

$$E|T_n|^p \leq (ET_n^2)^{p/2} = (\sum_{1\leq i_1<...<i_m\leq n} EY_{i_1i_2,...,i_m}^2)^{p/2} \leq n^{mp/2}(EY_{1,2,...,m}^2)^{p/2}. \qquad (19)$$

Переходя при оценке $E|W_{2k-2}|^{p/2}$ к условным математическим ожиданиям и используя тот факт, что на σ-алгебре $\sigma(X_{i_{k+1}},...,X_{i_m})$ $W_{2k-2}$ можно представить в виде сумм статистик порядка не превосходящего $2k-2$, ядро которых удовлетворяет условиям

(1)-(3) с заменой $p$ на $p/2$ и $E(\cdot)$ на $E(\cdot|X_{i_{k+1}},...,X_{i_m})$, из соотношений (16)-(19) легко получаем, что правое неравенство (5) верно при $2\leq p<4$.,,,,

Пусть $s(p)=[p/2]$. Предположим, что справедливость правой части (5) для всех симметрических статистик (любого порядка), удовлетворяющих условиям (1)-(3), и всех $p$, для которых $s(p)=l\geq 1$, уже доказана. Применяя (5) к условному математическому ожиданию относительно σ-алгебры $\sigma(X_{i_{k+1}},...,X_{i_m})$ каждой статистики, входящей в $W_{2k-2}$, $k=2,...,m$, используя неравенство Шварца и соотношения (7) для оценки полученных выражений по аналогии с [2], из (16)-(18) нетрудно получить, что

$$E|T_n|^p \leq B(p,m) \max_{k=0,m} n^{p/2(m-k)+k} E(E(Y^2(X_1,X_2,...,X_m)/X_1,...,X_k))^{p/2}$$

для всех $p$ с $s(p)=l+1$ и всех $Y \in G(p)$.

С учетом справедливости правой части (5) при $2\leq p<4$, то есть при $s(p)=1$, утверждение теоремы полностью доказано.

<u>Доказательство теоремы 3.</u> Из теоремы 1 следует, что

$$A(p,m) \max_{k=0,m} n^{p/2(m-k)+k} E(E(Y^2(X_1,X_2,...,X_m)/X_1,...,X_k))^{p/2} \leq$$

$$\leq E (\sum_{1\leq i_1<...<i_m\leq n} Y^2(X_{i_1},X_{i_2},...,X_{i_m}))^{p/2} \leq$$

$$\leq B(p,m) \max_{k=0,m} n^{p/2(m-k)+k} E(E(Y^2(X_1,X_2,...,X_m)/X_1,...,X_k))^{p/2}.$$

Используя неравенства (5), получаем (6).

R. Ibragimov, Sh. Sharahmetov

# MOMENT INEQUALITIES
# FOR SYMMETRIC STATISTICS

In this paper we prove the analoges of Khintchine and Rosenthal's moment inequalities for symmetric statistics of arbitrary order.

Р. Ибрагимов, Ш. Шарахметов

# МОМЕНТНЫЕ НЕРАВЕНСТВА
# ДЛЯ СИММЕТРИЧЕСКИХ СТАТИСТИК

В работе доказаны аналоги моментных неравенств Хинчина и Розенталя для симметрических статистик произвольного порядка.